\documentclass[12pt]{article}  
\usepackage{amssymb}
\usepackage{amsmath,amsthm}
\newtheorem{thm}{Theorem}[section]

\newtheorem{lem}[thm]{Lemma}

\newtheorem{prop}[thm]{Proposition}

\newtheorem{prob}[thm]{Problem}
\newtheorem{definitionS}[thm]{Definition}
\newenvironment{defn}{\begin{definitionS}\rm}{\end{definitionS}}
\newtheorem{example}[thm]{Example}

\newtheorem{remark}[thm]{Remark}
\newenvironment{rem}{\begin{remark}\rm}{\end{remark}}

\newcommand{\eqr}[1]{~\mbox{$(${\rm \ref{#1}}$)$}}
\newcommand{\Section}[1]{\section{#1}\setcounter{equation}{0}}

\newcommand{\ver}{\vspace*{-1.5mm}}

\newcommand{\cLG}{{\mathbb LG}}
\newcommand{\cZ}{{\mathcal Z}}
\newcommand{\C}{{\mathbb C}}
\newcommand{\K}{{\mathbb K}}
\newcommand{\R}{{\mathbb R}}
\newcommand{\pn}[1]{{\mathbb P}^{#1}}
\newcommand{\sym}{\mathrm{Sym}(n)}
\newcommand{\row}{{\rm rowsp}\,}
\newcommand{\tr}{\,{\rm trace}}

\newcommand{\G}{{\rm Grass}}

\newcommand{\vier}[4]{\left[ \begin{array}{ccc}
                   #1 &\;& #2 \\ #3 &\;& #4 \end{array} \right]}

\title{Output Feedback Pole Assignment for\\
  Transfer Functions with Symmetries\thanks{This work was carried
    out during the Special Year in Control 2003 at the
    Mittag-Leffler Institute, Stockholm, whose support is
    gratefully acknowledged.}}
\author{Uwe Helmke\\
  {\normalsize Department of Mathematics}\ver\\
  {\normalsize University of W\"urzburg}\ver\\
  {\normalsize Am Hubland}\ver\\
  {\normalsize 97074 W\"urzburg, Germany}\ver\\
  {\normalsize helmke@mathematik.uni-wuerzburg.de}\ver 
   \and 
  Joachim Rosenthal\thanks{Supported
    in part by NSF grant DMS-00-72383.}\\
{\normalsize Department of Mathematics}\ver\\
  {\normalsize University of Z\"urich}\ver\\
  {\normalsize Winterthurerstr 190}\ver\\
  {\normalsize CH-8057 Z\"urich,  Switzerland}\ver\\
  {\normalsize  http://www.math.unizh.ch/aa/}
  \and
  Xiaochang Wang\\
  {\normalsize Department of Mathematics and Statistics}\ver\\
  {\normalsize Texas Tech University}\ver\\
  {\normalsize Lubbock, TX 79409-1042}\ver\\
  {\normalsize alex.wang@ttu.edu}}

\begin{document}\maketitle

\begin{abstract}
  This paper studies the problem of pole assignment for
  symmetric and Hamiltonian transfer functions. A necessary and
  sufficient condition for pole assignment by complex symmetric
  output feedback transformations is given. Moreover, in the case
  where the McMillan degree coincides with the number of
  parameters appearing in the symmetric feedback transformations,
  we derive an explicit combinatorial formula for the number of
  pole assigning symmetric feedback gains. The proof uses
  intersection theory in projective space as well as a formula
  for the degree of the complex Lagrangian Grassmann manifold.
\end{abstract}
\medskip


\newpage
\Section{Introduction}

One of the best known inverse eigenvalue problems from linear
system theory is that of pole assignment, i.e. to find a static
output feedback gain for a given linear system such that the
closed loop poles of the system coincide with a specified subset
of the complex plane. Moreover, in the case of finitely many
solutions, a formula for the number of pole assigning feedback
transformations is desirable.  Early contributions on the subject
were obtained by e.g. Davison and Wang~\cite{da75} and
Kimura~\cite{ki75}, who derived sufficient conditions for the
solvability. However these conditons were far from being
necessary as well. In a series of pioneering
papers~\cite{he77,ma77a,ma78}, R. Hermann and C. F. Martin
applied tools from algebraic geometry to obtain necessary and
sufficient conditions, valid for a generic class of systems and
for {\em complex} feedback transformations. Their approach was
based on the dominant morphism theorem [Chapter AG, \S 17,
Theorem 17.3]~\cite{bo91b1} from complex algebraic geometry.  A
second breakthrough was subsequently made by R. W. Brockett and
C. I.  Byrnes~\cite{br81}, who used intersection theoretic
arguments and the Schubert calculus on Grassmann manifolds to
count the number of pole assigning complex feedback
transformations.  By refining these algebraic--geometric
approaches of Hermann and Martin, and Brockett and Byrnes, a
number of fundamental contributions on the subject were made that
finally led to a solution of the problem in the {\em real} case,
with important contributions due to~\cite{er02,le95,ro95,wa92}.
For an excellent survey paper on this subject, written from a
control-theoretic point of view, see e.g.~\cite{by89}. More
recently, various intersection theoretic tasks related to the
Schubert calculus have been studied in the algebraic geometry
literature; see e.g.~\cite{fu00,hu98,so00}.  The focus of most of
the investigations has been so far on the unstructured case,
where no underlying symmetries for the involved transfer function
or for the associated feedback transformations are imposed.
However, transfer functions with symmetries occur naturally in
various application areas, such as in network theory or
mechanics. For example, the transfer functions $G(s)$ of linear
$RLC$ - circuits, consisting solely of restistors, capacitors and
inductive elements are symmetric, i.e.  they satisfy
$G(s)^t=G(s)$. In mechanics, the transfer functions of linear
Hamiltonian systems are characterized by the symmetry relation
$G(-s)^t=G(s)$, while second order mechanical systems of the form
$$
M\ddot{x} = Nx +Bu, y=B^t x
$$
yield symmetric Hamiltonian transfer functions, satisfying
$$
G(s)=H(s^2), H(s)=H(s)^t;
$$
see e.g.~\cite{bi77a,by81,by80a1,fu83}.  For such structured
systems it is reasonable to restrict the class of admissible
feedback transformations to those that preserve the symmetry
properties of the transfer functions. Therefore the known results
on pole placement on unstructured systems do not apply in these
cases and require instead a new approach.

In this paper we start an investigation of the pole placement
problem for $n\times n$ symmetric transfer functions
$G(s)=G(s)^t$, arising in electrical network theory, and
Hamiltonian transfer functions. For both types of systems the
natural class of admissible output feedback tranformations are
the symmetric ones $F=F^t$, yielding a symmetric closed loop
transfer function
$$
G_F(s):=(I_n-G(s)F)^{-1}G(s).
$$
As the number of free parameters occuring in the symmetric
feedback matrices $F$ is $n(n+1)/2$, a necessary condition for
generic solvability of this output feedback problem is that the
McMillan degree $\delta$ of the transfer function $G$ satisfies
$\delta\geq \binom{n+1}{2}$ in the symmetric case, and
$\delta\geq n(n+1)$ in the Hamiltonian case. In fact, we show
that generically for {\em complex} symmetric output feedback
transformations this condition is also sufficient. Moreover, for
the limit case $\delta= \binom{n+1}{2}$ (or $\delta= n(n+1)$), we
derive an explicit combinatorial formula for the number of
complex symmetric output feedback gains that place the poles at
given points. Our formula coincides with that of the degree for
the complex Lagrangian manifold, given in~\cite{to03}.

In the real case such complete results can not be expected. In
fact, the symmetry of the transfer functions then imposes a
priori limitations for the possible pole locations of such
systems. This has been observed in~\cite{ma98}, where it is shown
for symmetric transfer functions that -- in the special case that
the Cauchy index of $G$ coincides with the McMillan degree --
then generically real symmetric output feedback pole
assignability holds if and only $n\geq \delta$. Of course, in
most applications we have $n\leq \delta$ and therefore the
description of the set of poles that can be achieved by real
symmetric output feedback becomes a complicated and nontrivial
task.
\Section{Complex symmetric and Hamiltonian realizations} \label{Sec-symm}

In this section we recall some basic facts concerning complex
symmetric and Hamiltonian transfer functions, respectively and
associated signature symmetric and Hamiltonian realizations. Let
$\C$ denote the field of complex numbers. A complex rational
transfer function $G(s) \in \C(s)^{n\times n}$ of McMillan degree
$\delta$ is called {\em symmetric, or Hamiltonian, respectively,}
if
$$
G(s)=G(s)^t, \quad \mathrm{or} \quad G(s) = G(-s)^t, \quad
\mathrm{respectively,}
$$
holds for all $s \in \C$.  A {\em complex
  symmetric realization} is a linear system of the form
$$
\dot{x}=Ax+Bu, y=B^t x,
$$
where $A\in \C^{\delta \times
  \delta}$ is symmetric, i.e. $A^t=A$, and $B\in \C^{\delta
  \times n}$. Similarly, a {\em Hamiltonian realization} is a
linear system
$$
\dot{x}=Ax+Bu, y=Cx,
$$
where $A\in \C^{\delta \times \delta},
B\in \C^{\delta \times n}, C\in \C^{n \times \delta}$ satisfies
$$
AJ = (AJ)^t, \quad C^t=JB
$$
and
$$
\vier{0}{I}{-I}{0}
$$
denotes the standard symplectic form on $\C^{\delta \times
  \delta}$.  In particular, Hamiltonian systems have always even
McMillan degree $\delta$.

Complex symmetric realizations are the natural class of
realizations for complex symmetric transfer functions. In fact,
they are the proper analogue of signature symmetric realizations
of real rational transfer functions, appearing in network theory.
Over $\R$, real symmetric realizations correspond to linear
models of $RC-$ networks, constructed entirely using capacitors
and resistors. The real symmetric transfer functions defined by
them are characterized by the property that the Cauchy-Maslov
index coincides with the Mcmillan
degree,~\cite{bi77a,fu83}.

The following variant of the Kalman realization theorem is
well-known; see e.g.~\cite{fu83,fu95}. Recall that the
complex orthogonal group $O(\delta,\C)$ is the matrix group
consisting of all complex $\delta \times \delta$ matrices $S$,
satisfying $SS^t=I_\delta$.  Given any complex realization
$(A,B,C)$ of a symmetric transfer function $G(s)=C(sI-A)^{-1}B$,
note that $(A^t,C^t,B^t)$ is also a realization.

\begin{prop}
  Let $G(s)=G(s)^t$ be an $n \times n$ strictly proper, complex
  rational transfer function of McMillan degree $\delta$. Then
\begin{itemize}
\item[(1)] $G(s)$ has a controllable and observable complex
  symmetric realization $(A,B,C)=(A^t,C^t,B^t)$.
\item[(2)] If $(A_i,B_i,C_i)$, i=1,2, are two controllable and
  observable complex symmetric realizations of $G(s)$, then there
  exists a unique complex orthogonal transformations $S\in
  O(\delta,\C)$ such that $(A_2,B_2,C_2)=(SA_1 S^{-1},SB_1,C_1
  S^{-1})$.
\end{itemize}
\end{prop}
In the literature usually only the real case of the above result
is proven, where the statement is actually slightly different due
to the presence of signature symmetric realizations. In the
complex case the result simplifies to the one given here.  For
the sake of completeness we include the proof; see also
\cite{fu95}.
\begin{proof}
  If $(A,B,C)$ is a minimal realization of $G(s)$ then, by
  symmetry of $G$, also $(A^t,C^t,B^t)$ is a minimal realization.
  Applying Kalman's realization theorem implies the existence of
  a unique invertible complex $\delta \times \delta$ matrice $S$
  with
  $$
  (A^t,C^t,B^t)=(SAS^{-1},SB,CS^{-1}).
  $$
  By transposing this
  equation and using the uniqueness of $S$ we conclude that
  $S=S^t$. It is a well known fact from linear algebra that every
  complex symmetric invertible matrix has a representation
  $S=XX^t$ by a complex invertible matrix $X$. Moreover, $X$ is
  uniquely determined up to right factors $XT$ where $T\in
  O(\delta,\C)$. Then $(XAX^{-1},XB,CX^{-1})$ is a complex
  symmetric realization, which completes the proof.
\end{proof}
There is a similar realization theorem for Hamiltonian systems,
for which we refer to the literature; see
e.g.~\cite{by80a1,fu83}.
Static linear output feedback can be meaningfully defined for
such systems only through symmetric gain matrices.  Thus an
output feedback transformation
$$
u=Fy+v
$$
with the closed system
$$
\dot{x}=(A+BFB^t)x+Bu, y=B^t x
$$
preserves the complex symmetry
of the realizations if and only if $F=F^t$. Thus we define two
complex symmetric realizations $(A_i,B_i,C_i)$ to be {\em
  symmetric output feedback equivalent} if and only if there
exist $S\in O(\delta, \C)$, $F=F^t \in \C^{n \times n}$ with
$$
(A_2,B_2,C_2)=(S(A_1+B_1FB_1^t)S^{-1},SB_1,C_1S^{-1}).
$$
Equivalently, if and only if for the associated transfer
functions $G_i(s)$:
$$
G_2(s):=(I_n-G_1(s)F)^{-1}G_1(s).
$$
Similarly, output feedback for Hamiltonian systems
$$
\dot{x}=(A+BFC)x+Bu, y=Cx
$$
preserves the Hamiltonian
properties of the realization if and only if $F=F^t$. Thus in
both cases we have to focus on symmetric output feedback.

We note some elementary geometric properties of the set of
complex symmetric transfer functions that will be important in
the subsequent development; see e.g. \cite{by80a1} for providing some
of the details for the proof of the subsequent theorem. We omit a full
proof as it would take us to far apart from the subject.

\begin{prop}
  Let $SRat_{\delta,n}(\C)$ and $Ham_{\delta,n}(\C)$,
  respectively denote the sets of strictly proper, complex
  symmetric and Hamiltonian, respectively, $n\times n$ transfer
  functions of McMillan degree $\delta$. Then
  $SRat_{\delta,n}(\C)$ and $Ham_{\delta,n}(\C)$, respectively,
  is a smooth complex manifold of complex dimension
  $\delta(n+1)$, and dimension $\delta n$ respectively. Moreover,
  they are nonsingular irreducible quasi-affine varieties.
\end{prop}
In particular, there is a canonical notion of ``genericity'' for
complex symmetric or Hamiltonian transfer functions. Explicitely,
a property $E$ of complex symmetric transfer functions is called
generic, if the set defined by $E$
$$
\{G\in SRat_{\delta ,n}(\C)\; | \; G\; has\; property\; E \}
$$
is a Zariski-open subset of $SRat_{\delta,m}(\C)$. Equivalently,
this can be also expressed in terms of complex symmetric
realizations.

\Section{Main result} \label{Sec-compact}

After these preliminaries we can now rigorously formulate and
proof the main technical results of this paper.  Let $G(s)$ be an
$n\times n$ complex symmetric or Hamiltonian transfer function,
i.e. $G(s)^t=G(s)$ or $G(-s)^t=G(s)$, respectively. Assume that
$G(s)$ is strictly proper and has McMillan degree $\delta$. The
complex symmetric eigenvalue assignment problem then asks the
following question:
\begin{prob}                        \label{problem}
  Given an arbitrary monic polynomial $\varphi(s)\in\C[s]$ of
  degree $\delta$ ($\varphi(s)=\varphi(-s)$ is assumed to be {\bf
    even} in the Hamiltonian case). Is there an $n\times n$
  complex symmetric matrix $F$ such that the closed loop transfer
  function
  $$
  G_F(s):=(I_n-G(s)F)^{-1}G(s)
  $$
  has characteristic polynomial $\varphi(s)$, i.e. the poles
  of $G_F(s)$ are the zeroes of $\varphi(s)$?
\end{prob}
If for a particular symmetric (Hamiltonian) transfer function
$G(s)$ Problem~\ref{problem} has an affirmative answer we will
say that $G(s)$ is {\em pole assignable in the class of complex
  symmetric (Hamiltonian) feedback compensators}. We say that
$G(s)$ is {\em generically pole assignable}, if the problem is
solvable for a generic choice of admissible polynomials
$\varphi(s)$.

Similar to the situation of the static pole placement
problem~\cite{br81,wa92} and the dynamic pole placement
problem~\cite{ro94}, Problem~\ref{problem} turns out to be highly
nonlinear and techniques from algebraic geometry will be required
to study the problem.  The first main result is in the spirit of
Hermann and Martin, by deriving a generic necessary and
sufficient condition via the dominant morphism theorem.

We prove some lemmas first. Let
$\pi(A)=(a_{11},\dots,a_{\delta\delta})$ be the projection onto
the diagonal entries of an $\delta\times \delta$ matrix $A$. In
the sequel we will identify $\mathbb{C}^\delta$ with the complex
vector space of row vectors. For any symmetric matrix $L$, define
$\theta^{L}: O(\delta, \mathbb{C})\rightarrow \mathbb{C}^\delta$
through
$$
\theta^{L} (S)= \pi(SLS^{-1}).
$$
As $O(\delta, \mathbb{C})$ is a Lie group, its tangent space
at the identity matrix $I$ is given by the Lie algebra of complex
skew-symmetric matrices
$$
so(\delta, \mathbb{C})=\{X\in \mathbb{C}^{\delta\times \delta}
\mid X+X^t=0\}.
$$
Moreover, the Jacobian $d\theta^{L}_{I}$ of $\theta^{L}$ at
$I$ is given by
$$
d\theta^{L}_{I}: so(\delta,\mathbb{C}) \rightarrow V,\ \ \ \ \ 
d\theta^{L}_{I}(X)= \pi(XL-LX),
$$
where
$$
V=\{(x_1,\dots,x_\delta) \in \mathbb{C}^{\delta} \mid
\sum_{1}^\delta x_i=0\}.
$$

For any $\delta\times \delta$ matrix $L$, the graph ${\cal G}(L)$
of $L$ is defined as a graph with $\delta$ vertices such that
there is a path from vertex $i$ to vertex $j$ if and only if the
$ij$th entry of $L$ is none zero. It is a well-known fact from
linear algebra, that the graph ${\cal G}(L)$ is connected if and
only if $L$ is irreducible, i.e. if and only if there exists no
permutation matrix $P$ such that $PLP^{-1}$ is block diagonal. We
use this fact together with an idea developed
in~\cite[Lemma~2.5]{he97} to prove the following equivalent
characterization.

\begin{lem}
  The Jacobian $d\theta^{L}_{I}$ is surjective if and only if the
  associated graph ${\cal G}(L)$ is connected.
\end{lem}

\begin{proof}
  By inspection, the derivative $d\theta^{L}_{I}$ is not
  surjective if and only if there exists a nonzero diagonal
  matrix $Z$ of trace zero, such that for all $X\in so(\delta,
  \mathbb{C})$
  $$
  \tr (Z(XL-LX))=\tr ((LZ-ZL)X)=0.
  $$
  By symmetry of $L,Z$ we have $LZ-ZL\in so(\delta,
  \mathbb{C})$.  Since the trace function defines a nondegenerate
  bilinear form on $so(\delta, \mathbb{C})$, the condition $\tr
  ((LZ-ZL)X)=0$ is equivalent to $LZ=ZL$. Since $Z$ is a nonzero
  diagonal matrix of trace zero, there is a permutation matrix
  $P$ such that $\hat{Z}:=PZP^{-1}={\rm block\ 
    diag}(a_1I_1,\dots,a_kI_k)$ with $k\geq 2$ and $a_i$'s
  distinct. Let $\hat{L}=PLP^{-1}$. Then $LZ=ZL$ is equivalent to
  $\hat{L}\hat{Z}=\hat{Z}\hat{L}$, which is equivalent to
  $\hat{L}$ being block diagonal. But from the above remark this
  is equivalent to the graph ${\cal G}(L)$ being disconnected.
  The result follows.
\end{proof}

\begin{lem} \label{lem1}
  Let $L$ be a nonzero complex symmetric matrix such that
  $\pi(L)=0$.  Then there is a family of orthogonal matrices
  $S(\epsilon)\in O(\delta, \mathbb{C})$, $\epsilon\geq 0$, with
  $S(0)=I$ such that the matrix
  $\hat{L}(\epsilon):=S(\epsilon)LS(\epsilon)^{-1}$ has the
  properties that $\pi(\hat{L}(\epsilon))=0$ and
  $d\varphi^{\hat{L}(\epsilon)}_{I}$ is surjective for all
  $\epsilon\in (0,\pi/2)$.
\end{lem}

\begin{proof}
  If ${\cal G}(L)$ is connected, then by the previous lemma the
  choice $S(\epsilon):=I$ does the job. Thus it suffices to prove
  that ${\cal G}(L)$ not connected implies that then we can find
  a family of transformations $S(\epsilon)$, such that
  $\pi(S(\epsilon)LS(\epsilon)^{-1})=0$ and the largest connected
  subgraph of ${\cal G}(S(\epsilon)LS(\epsilon)^{-1})$ contains
  more vertices than that of ${\cal G}(L)$ for all
  $0<\epsilon<\pi/2$.
  
  Note that $\pi(L)=0$ and $L\neq 0$ imply that the largest
  connected subgraph of ${\cal G}(L)$ must contain at least $2$
  vertices.  Assume that the largest connected subgraph of ${\cal
    G}(L)$ contains $k$ vertices, $2\leq k <\delta$.  Without
  loss of generality, assume
  $$
  L=\left[\begin{array}{cc}
      L_1& 0\\
      0 &L_2
\end{array}\right]
$$
where the graph of the $k\times k$ sub-matrix $L_1$ is
connected.  Write
$$
L_1=\left[\begin{array}{cc}
    L_{11}& \alpha\\
    \alpha^t & 0
\end{array}\right]\mbox{ and }L_2=\left[\begin{array}{cc}
0 & \beta^t\\
\beta & L_{22}
\end{array}\right],$$
where $L_{11}$ and $L_{22}$ are sizes $(k-1)\times (k-1)$ and
$(\delta-k-1)\times (\delta-k-1)$ respectively. By irreducibility
of $L_1$ we have $\alpha\neq 0$.  Thus $L$ has the form
$$\left[\begin{array}{cccc}
    L_{11}   & \alpha & 0       & 0    \\
    \alpha^t & 0      & 0       & 0    \\
    0        & 0      & 0       &\beta^t \\
    0 & 0 & \beta &L_{22}
\end{array}\right].
$$
Let
$$
S(\epsilon)=\left[\begin{array}{cccc}
    I_{(k-1)\times (k-1)}   & 0      & 0  & 0    \\
    0   & \cos \epsilon & -\sin  \epsilon & 0    \\
    0   & \sin \epsilon & \cos \epsilon   & 0    \\
    0 & 0 & 0 & I_{(\delta-k-1)\times (\delta-k-1)}
\end{array}\right].$$
Then
$$
S(\epsilon)LS^{-1}(\epsilon)=\left[\begin{array}{cccc}
    L_{11}                  & (\cos \epsilon)\alpha  & (\sin \epsilon)\alpha      &  0                      \\
    (\cos \epsilon)\alpha^t & 0                      & 0                          & (-\sin \epsilon)\beta^t   \\
    (\sin \epsilon)\alpha^t & 0                      & 0                          & (\cos \epsilon)\beta^t    \\
    0 &(-\sin \epsilon)\beta &(\cos \epsilon)\beta & L_{22}
\end{array}\right].
$$
For the graph ${\cal G}(S(\epsilon)LS^{-1}(\epsilon))$ with
$0<\epsilon < \pi/2$, the vertices $\{1, \dots, k\}$ are still
connected and the vertex $k+1$ is symmetrically connected to at
least one of the first $k$ vertices. Thus the vertices $\{1,
\dots, k+1\}$ are connected.
\end{proof}

\begin{lem}\label{lem3}
  Let ${\cal L}$ be a linear subspace of complex symmetric
  matrices of dimension $\delta$, and ${\cal L}\not\subset
  sl(\delta,\mathbb{C})$.  Then there exists an orthogonal matrix
  $S\in O(\delta,\mathbb{C})$ such that $\pi\mid_{S{\cal
      L}S{-1}}$ is one to one, and onto.
\end{lem}

\begin{proof}
  The proof goes by recursively constructing a basis
  $\{L_1,\dots,L_\delta\}$
  of ${\cal L}$ such that\\
  $\{\pi(SL_1S^{-1}),\dots,\pi(SL_\delta S^{-1})\}$ are linearly
  independent for a suitable complex orthogonal matrix $S\in
  O(\delta,\mathbb{C})$.  First note, that we can modify any
  basis of ${\cal L}$ into a basis ${\cal
    L}^{(1)}:=\{L_1,\dots,L_\delta\}$ of ${\cal L}$ such that
  $L_1\not\in sl(\delta,\mathbb{C})$, and $L_i\in
  sl(\delta,\mathbb{C})$, for $i=2,\dots,\delta$. In fact, if
  $\{K_1,\dots,K_\delta\}$ denotes any basis of ${\cal L}$ with
  $\tr ( K_1)\neq 0$, then
  $\{L_1:=K_1,L_2:=K_2-c_2K_1,\dots,L_\delta:=K_\delta-c_\delta
  K_1\}$, $c_i:=\tr (K_i)/\tr (K_1)$, is as desired. By
  construction of $L_1$, then $\pi(L_1)\neq 0$.
  
  Let $\{L_1,\dots,L_\delta\}$ be a basis of ${\cal L}$ such that
  $L_1\not\in sl(\delta,\mathbb{C})$, and $L_i\in
  sl(\delta,\mathbb{C})$, for $i=2,\dots,\delta$.  Then ${\rm
    dim\ span} \{\pi(L_1),\dots,\pi(L_\delta)\}:=k \geq 1$.  If
  $k<\delta$, then by re-ordering the indices we can assume that
  $\{\pi(L_1),\dots,\pi(L_k)\}$ are linearly independent, and
  $$
  \pi(L_j)= \sum_{i=1}^{k}c_{ij}\pi(L_i)\mbox{ for
    $j=k+1,\dots,\delta$.}
  $$
  By replacing $L_j$ with $L_j-\sum_{i=1}^{k}c_{ij}L_i$ we can
  further assume that $\pi(L_j)=0$ for $j=k+1,\dots,\delta$.  It
  is thus sufficient to show that if there is an orthogonal
  matrix $\hat{S}$ such that the matrices
  $\{M_j:=\hat{S}L_j\hat{S}^{-1},j=1,\dots,\delta\}$ have the
  property that $\{\pi(M_1),\dots,\pi(M_k)\}$ are linearly
  independent, and $\pi(M_j)=0$, $j=k+1,\dots,\delta$, for some
  $k<n$, then we can find an orthogonal $S$ such that

  $$
  \{\pi(SM_1S^{-1}),\dots,\pi(SM_kS^{-1}),\pi(SM_{k+1}S^{-1})\}
  $$
  are linearly independent.
  
  By Lemma~\ref{lem1}, there exists $S_\epsilon \in
  O(\delta,\mathbb{C})$ arbitrarily close to the identity matrix
  such that
  $\pi(S_{\epsilon}M_1S_{\epsilon}^{-1}),\dots,\pi(S_{\epsilon}M_kS_{\epsilon}^{-1})$
  are linearly independent and the graph ${\cal
    G}(S_{\epsilon}M_{k+1}S_{\epsilon}^{-1})$ is connected. By
  replacing $M_i$ with $S_{\epsilon}M_iS_{\epsilon}^{-1}$ we can
  assume further that $d\theta^{M_{k+1}}_{I}$ is onto $V$.  Then
  there exists a skew-symmetric matrix $X$ such that
  $$
  \pi(XM_{k+1}-M_{k+1}X)\not\in {\rm
    span}\{\pi(M_1),\dots,\pi(M_k)\}.
  $$
  Let
  $$
  S(\epsilon)=\exp (\epsilon X).
  $$
  Then $S(\epsilon)$ is orthogonal for all $\epsilon$, and
  $$
  S(\epsilon)=I+\epsilon X+\mbox{higher order terms}.
  $$
  The Taylor series expansions of $\{\pi(S(\epsilon)M_i
  S(\epsilon)^{-1})\}$ have the forms
  $$
  \pi(S(\epsilon)M_i
  S(\epsilon)^{-1})=\pi(M_i)+\beta_i(\epsilon),\ \ \ i=1,\dots,k,
  $$
  and
  $$
  \pi(S(\epsilon)M_{k+1} S(\epsilon)^{-1})=\epsilon
  \left(\pi(XM_{k+1}-M_{k+1}X)+\beta_{k+1}(\epsilon)\right)
  $$
  where $\beta_i(\epsilon)$ are continuous with respect to
  $\epsilon$ and $\beta_i(\epsilon)\rightarrow 0$ as
  $\epsilon\rightarrow 0$ for $i=1,\dots,k+1$.  Since
  $\{\pi(M_1),\dots,\pi(M_k),\pi(XM_{k+1}-M_{k+1}X)\}$ are
  linearly independent, for sufficient small $\epsilon>0$,
  $\{\pi(M_1)+\beta_1(\epsilon),\dots,\pi(M_k)+\beta_k(\epsilon),
  \pi(XM_{k+1}-M_{k+1}X)+\beta_{k+1}(\epsilon)\}$ are also
  linearly independent, i.e.
  $\{\pi(S(\epsilon)M_1S(\epsilon)^{-1},\dots,
  \pi(S(\epsilon)M_kS(\epsilon)^{-1},\pi(S(\epsilon)M_{k+1}S(\epsilon)^{-1}\}$
  are linearly independent.
\end{proof}

\begin{thm}                                     \label{premain}
  If $G(s)$ is a symmetric (or Hamiltonian) transfer function of
  McMillan degree $\delta> \binom{n+1}{2}$ (or $\delta> n(n+1)$),
  then $G(s)$ is not pole assignable in the class of (real or)
  complex symmetric feedback compensators.
  
  When $\delta\leq \binom{n+1}{2}$ (or $\delta\leq n(n+1)$), then
  there is a generic set of $n\times n$ symmetric (or
  Hamiltonian) transfer functions of degree $\delta$ which are
  generically pole assignable via complex symmetric feedback
  compensators.
\end{thm}

\begin{proof}
  We only give a sketch of the proof, as the arguments based on
  the dominant morphism theorem are well known from
  \cite{he77,ma77a}. Note, however, that there is serious gap in
  the proof of \cite{ma77a} for the pole placement result on
  Hamiltonian systems because it is not proved that the set of
  generically pole assignable Hamiltonian systems is non empty.
  In fact, a construction of such an example is not completely
  trivial and depends on our previous lemmas.
  
  The first claim follows immediately from a standard dimension
  argument, as the vector space $\sym$ of complex $n \times n$
  symmetric matrices has dimension $\binom{n+1}{2}$.  For the
  second claim we note that the set of generically pole
  assignable systems is a Zariski open subset of the nonsingular,
  irreducible quasi-affine variety of symmetric or Hamiltonian
  transfer functions, respectively. Therefore we only need to
  show that this Zariski open subset is nonempty. By the Dominant
  Morphism Theorem, it suffices to find one system whose Jacobian
  of the pole placement map at one point is onto.
  
  Note, by the Newton formula, that the coefficients of the
  characteristic polynomial $\det
  (sI-A)=s^\delta+\alpha_{\delta-1}s^{\delta-1}+ \cdots \alpha_1
  s+\alpha_0$ are related to the traces of powers of $A$ as
  follows:
\begin{eqnarray*}
\alpha_{\delta-1}&=&-\tr(A)\\
\alpha_{\delta-2}&=&-\frac{1}{2}(\tr(A^2)+\alpha_{\delta-1} \tr(A))\\
&\vdots&\\
\alpha_0&=&-\frac{1}{\delta}(\tr(A^\delta)+\alpha_{\delta-1}
\tr(A^{\delta-1})+\cdots+\alpha_1 \tr(A)).
\end{eqnarray*}

Therefore for the case of symmetric transfer functions, the pole
placement map is equivalent to the map
\begin{equation}                       \label{poleplace}
   \begin{array}{cccc}
   \phi: &\sym& \longrightarrow &{\mathbb C}^\delta \\
    & F & \longmapsto &(\tr(A+BFB^t), \cdots,\tr(A+BFB^t)^{\delta})
   \end{array}
\end{equation}
and its Jacobian at $0$ is given by
$$
d\phi_0 (F)=(\tr(BFB^t),2\tr(ABFB^t),\dots,\delta
\tr(A^{\delta-1}BFB^t).
$$

For the case of Hamiltonian transfer functions, since $JAJ=A^t$
and $J^2=-I$, we have $(-1)^{k-1}JA^kJ=(A^k)^t$ for
$k=1,2,\dots$, which implies that the characteristic polynomial
of $A$ is even and
$$
\tr (A^k)=0 \mbox{ holds for all odd $k$'s.}
$$
Therefore the pole placement map is equivalent to the map
\begin{equation}                       \label{poleplace2}
   \begin{array}{cccc}
   \psi: &\sym& \longrightarrow &{\mathbb C}^{\delta/2} \\
    & F & \longmapsto &(\tr(A+BFC)^2, \tr(A+BFC)^4,\cdots,\tr(A+BFC)^{\delta})
   \end{array}
\end{equation}
and its Jacobian at $0$ is given by
$$
d\psi_0 (F)=(2\tr(ABFC),4\tr(A^{3}BFC)\dots,\delta
\tr(A^{\delta-1}BFC).
$$

We first consider the case of symmetric transfer functions.  Let
$B$ be any real nonzero matrix and ${\cal L}=\{BFB^t\mid F\in
\sym\}$.  Then ${\cal L}\not\subset sl(\delta,\mathbb{C})$ and
${\rm dim} {\cal L}\geq \delta$. By Lemma~\ref{lem3} there exists
an orthogonal matrix $S\in O(\delta,\mathbb{C})$ such that
$\pi\!\!\mid_{S{\cal L}S^{-1}}$ is surjective. Let $D={\rm
  diag}(1,2,\dots,\delta)$ and $A=S^{-1}DS$.  Then
$$
\begin{array}{lcl}
d\phi_0 (F)&=&(\tr(BFB^t),2\tr(ABFB^t),\dots,\delta 
\tr(A^{\delta-1}BFB^t))\\
&=&(\tr(SBFB^tS^{-1}),2\tr(DSBFB^tS^{-1}),\dots,\delta
\tr(D^{\delta-1}SBFB^tS^{-1}))\\
&=&\pi(SBFB^tS^{-1})V
\end{array}
$$
where
$$
V=\left[\begin{array}{cccc}
    1&1&\cdots&1\\
    1&2&\cdots&2^{\delta-1}\\
    \vdots&\vdots&&\vdots\\
    1&\delta&\cdots&\delta^{\delta-1}
\end{array}\right]\left[\begin{array}{cccc}
1&0&\cdots&0\\
0&2&\cdots&0\\
\vdots&\vdots&&\vdots\\
0&0&\cdots&\delta
\end{array}\right].
$$
Since $\pi\!\!\mid_{S{\cal L}S^{-1}}$ is surjective and $V$ is
nonsingular, $d\phi_0$ is onto.

For the case of Hamiltonian transfer functions, Let
$$
B=\left[\begin{array}{c}
    0\\
    B_1
\end{array}\right]\mbox{ and } C=\left[\begin{array}{cc}
B^t_1& 0
\end{array}\right]
$$
where $B_1$ is any real nonzero $\frac{\delta}{2}\times n$
matrix, and ${\cal L}=\{B_1FB^t_1\mid F\in \sym\}$.  Then ${\cal
  L}\not\subset sl(\delta/2,\mathbb{C})$.  By Lemma~\ref{lem3}
there exists an orthogonal matrix $S_1\in O(\delta/2,\mathbb{C})$
such that $\pi: S_1{\cal L}S^{-1}_1\mapsto {\mathbb
  C}^\frac{\delta}{2}$ is surjective.  Let $D_1={\rm
  diag}(1,2,\dots,\delta/2)$,
$$
S=\left[\begin{array}{cc}
    S_1&0\\
    0&S_1
\end{array}\right], \mbox{ and } D=\left[\begin{array}{cc}
0&D_1\\
D_1&0
\end{array}\right],
$$
and $A=S^{-1}DS$. Note that $D,S$ are Hamiltonian and
symplectzic matrices, respectively. In particular, $A$ is
Hamiltonian.  Then
$$
\begin{array}{lcl}
d\psi_0 (F)&=&(2\tr(ABFC),4\tr(A^3BFC),\dots,\delta \tr(A^{\delta-1}BFC))\\
&=&(2\tr(DSBFCS^{-1}),4\tr(D^3SBFCS^{-1}),\dots,\delta
\tr(D^{\delta-1}SBFCS^{-1}))\\
&=&(2\tr(D_1S_1B_1FB^t_1S_1^{-1}),4\tr(D_1^3S_1B_1FB^t_1S_1^{-1}),\dots,\delta
\tr(D_1^{\delta-1}S_1B_1FB^t_1S_1^{-1}))\\
&=&\pi(S_1B_1FB^t_1 S_1^{-1})U
\end{array}
$$
where
$$
U=\left[\begin{array}{cccc}
    1&1&\cdots&1\\
    2&2^3&\cdots&2^{\delta-1}\\
    \vdots&\vdots&&\vdots\\
    \frac{\delta}{2}&\left(\frac{\delta}{2}\right)^3&\cdots&
    \left(\frac{\delta}{2}\right)^{\delta-1}
\end{array}\right]\left[\begin{array}{cccc}
2&0&\cdots&0\\
0&4&\cdots&0\\
\vdots&\vdots&&\vdots\\
0&0&\cdots&\delta
\end{array}\right].
$$
Since $\pi\!\!\mid_{S_1B_1FB^t_1S_1^{-1}}$ is surjective and
$U$ is nonsingular, $d\psi_0$ is onto.
\end{proof}

The second main theorem in this paper deals with the limit case
$\delta= \binom{n+1}{2}$, where we can prove a more precise
statement.

\begin{thm}                                     \label{main}
  Let $\delta= \binom{n+1}{2}$ in the symmetric case, and
  $\delta= n(n+1)$ for Hamiltonian systems. Then for a generic
  set of $n\times n$ symmetric (or Hamiltonian) transfer
  functions of degree $\delta$ the number of pole assigning
  complex symmetric feedback compensators is finite and when
  counted with multiplicities there are exactly
  \begin{equation}                      \label{deg-formula}
    d(n) := 2^{\binom{n}{2}}\frac{\binom{n+1}{2}!\ 1!\ 2!\ \cdots \
      (n-1)!}{1!\ 3!\ \cdots\  (2n-1)!}=
{\frac {\binom{n+1}{2} !}{\prod _{i=0}^{n-1}
 \left( 2\,i+1 \right) ^{n-i}}}
 \end{equation}
 many symmetric compensators as solution.
\end{thm}



One immediately computes $d(1) = 1$, $d(2) = 2$, $d(3) = 2^4$,
$d(4) = 3 \cdot 2^8$, $d(5) = 11 \cdot 13 \cdot 2^{11}$ and $d(6)
= 13 \cdot 17 \cdot 19 \cdot 2^{18}$. The integer sequence $d(n)$
is sequence A005118 in Sloane's data bank of integer
sequences~\cite{sl03}.  The sequence has several combinatorial
and geometric interpretations. For the context of this paper it
will be important that $d(n)$ is equal to the the degree of the
Lagrangian Grassmannian, the projective variety of all maximal
isotropic subspaces in a complex vector space of dimension $2n$
and this has been recently established by Totaro~\cite{to03}.

As it can be seen from this sequence, $d(n)$ appears always to be
even, except for $n=1$. This is related to the fact, that the
symmetric output feedback pole placement problem is {\em not}
generically solvable over the reals. Actually more is true. The
sequence
$$
\tilde{d}(n):=d(n)2^{-\binom{n}{2}}
$$
is the degree of the spinor variety, the complex projective
variety $SO(2n+1)/U(n)$~\cite{hi82a}; in particular
$\tilde{d}(n)$ represents an integer sequence again.  The
sequence $\tilde{d}(n)$ appears under the number A003121 in
Sloane's data bank~\cite{sl03}.\medskip

The proof of Theorem~\ref{main} will occupy the rest of this
section. The proof will necessitate a geometric reformulation and
several technical lemmas.

First we will describe the closed loop characteristic equation in
a slightly more convenient way. Consider a left coprime
factorization $D^{-1}(s)N(s)=G(s)$ of the symmetric or
Hamiltonian transfer function $G(s)$. Let $F\in\sym$ be an
$n\times n$ complex symmetric matrix.  When the feedback law
$y=-Fu+v$ is applied then up to a constant factor the
characteristic polynomial $\varphi(s)$ is also equal to
\begin{equation}                   \label{clo-sys}
     \det\vier{D(s)}{N(s)}{F}{I_n}.
\end{equation}

The vector space $\sym$ describing the set of $n\times n$ complex
symmetric matrices is not very well suited to invoke strong
theorems from algebraic geometry and intersection
theory~\cite{fu84}, as these usually require compactness
assumptions on the underlying spaces. A similar difficulty exists
for the static output pole placement problem.  Brockett and
Byrnes showed in~\cite{br81} how to translate the static pole
placement problem into a geometric problem. This then resulted
into an intersection problem on a compact Grassmann variety and
methods from classical Schubert calculus~\cite{sc1879,so01} could
be invoked.

We will follow this compactification strategy for
Problem~\ref{problem} as well.  This will lead us to an
intersection problem on some projective variety.  In order to do
so we therefore need a good compactification of $\sym$.  For this
identify the rowspan $\row [F\ I_n]$ of any symmetric matrix $F$
with an element of the Grassmann variety $\G(n,\C^{2n})$. Using
the Pl\"ucker embedding
$$
\G(n,\C^{2n}) \longrightarrow {\mathbb P}\left(\wedge^{n}
  \C^{2n}\right)=\pn{N},\ \ N=\binom{2n}{n}-1
$$
we can then identify $\sym$ with a quasi-projective subset of
the complex projective variety $\pn{N}$.
\begin{defn}
  The algebraic closure of the set
  $$
  \left\{ \row [F\ I_n] \mid F\in\sym\right\}
  $$
  is called the {\em complex Lagrangian Grassmann manifold}.
  It will be denoted by $\cLG (n)$.
\end{defn}

It is well known that $\cLG (n)$ is a smooth projective variety
of of dimension $\binom{n+1}{2}$, the dimension of $\sym$.  Note
that every element in $\cLG (n)$ can be simply represented by a
subspace of the form $\row [F_1\ F_2]$, where $F_1(F_2)^t$ is a
symmetric matrix, i.e. $F_1(F_2)^t=F_2(F_1)^t$. The elements of
$\cLG (n)$ are thus exactly the Lagrangian subspaces of
$\C^{2n}$.  The subspace $\row [F_1\ F_2]$ coincides with the
subspace $\row [S\ I_n]$ associated with an element $S$ of $\sym$
if and only if $F_2$ is invertible. Moreover, then
$S=(F_2)^{-1}F_1$.  When $F_2$ is singular one can still define a
characteristic polynomial through
\begin{equation}                   \label{clo-sys2}
   \varphi(s):=  \det\vier{D(s)}{N(s)}{F_1}{F_2}.
\end{equation}

Note that in the Hamiltonian case, $\varphi(s)$ is necessarily
even, i.e. then
$$
\varphi(s)= \varphi(-s).
$$
Let $f_i, \; i=0,\ldots,N$ be the Pl\"ucker coordinates of
$\row [F_1\ F_2]$.  In terms of the Pl\"{u}cker coordinates the
characteristic equation can then be written as:
\begin{equation}    \label{eqw2}
\det\vier{D(s)}{N(s)}{F_1}{F_2}=
\sum_{i=0}^N p_i(s)  f_i,
\end{equation}
where $p_i(s)$ is the cofactor of $f_i$ in the
determinant~(\ref{eqw2}).

Let $\cZ \subset \pn{N}$ be the linear subspace defined by
\begin{equation}    \label{center}
\cZ =\{z\in {\mathbb P}^N|\sum_{i=0}^N p_i(s)z_i=0\}.
\end{equation}
Following~\cite{ki04,ra96,ro94,wa92} we identify a closed loop
characteristic polynomial $\varphi(s)$ with a point in
$\pn{\delta}$.  In analogy to the situation of the static pole
placement problem considered in~\cite{br81,wa92} (compare also
with~\cite[Section~5]{ro94}) one has a well defined
characteristic map
\begin{equation}                       \label{central}
   \begin{array}{cccc}
   \chi: &\cLG (n)-\cZ & \longrightarrow &{\mathbb P}^\delta \\
    & \row [F_1\ F_2] & \longmapsto &\sum_{i=0}^N f_ip_i(s).
   \end{array}
\end{equation}
in the complex symmetric case and
\begin{equation}                       \label{centralnew}
   \begin{array}{cccc}
   \chi': &\cLG (n)-\cZ & \longrightarrow &{\mathbb P}^{\delta /2}\\
    & \row [F_1\ F_2] & \longmapsto &\mbox{ even part of }\sum_{i=0}^N f_ip_i(s).
   \end{array}
\end{equation}
in the Hamiltonian case. In the latter case the reduction in
dimension of the projective space arises due to the evenness of
the closed loop characteristic polynomial, so that in the second
map only the coefficients of the even terms of $\sum_{i=0}^N
f_ip_i(s)$ do appear.

Recall the notion of degree of a variety~\cite[Chapter I, \S
7]{ha77} and the notion of a central projection
(see~\cite[Chapter I, \S 4]{sh77}).  The geometric properties of
the map $\chi$ are as follows:
\begin{thm}                             \label{thm2.3}
  The maps $\chi,\chi'$ define central projections. In particular
  if $\cZ \cap \cLG (n)=\emptyset$ and $\dim\cLG
  (n)=\binom{n+1}{2}=\delta$ then $\chi$ is surjective, and there
  are $\deg \cLG (n)$ many pre-images (counted with multiplicity)
  for each point in $\pn{\delta}$, where $\deg \cLG (n)$ is the
  degree of the Lagrangian manifold $\cLG (n)$ in ${\mathbb
    P}^N$.  Similarly, if $\dim\cLG (n)=\binom{n+1}{2}=\delta/2$
  then $\chi'$ is surjective with exactly $\deg \cLG (n)$ many
  pre-image points in each fiber.
\end{thm}
\begin{proof}
  By definition (see e.g.~\cite{mu76,sh77}) $\chi$ represents a
  central projection of $\cLG (n)$ {\em from the center} $\cZ$
  {\em to} $\pn{\delta}$. When $\cZ \cap \cLG (n)=\emptyset$ and
  $\dim\cLG (n)=\binom{n+1}{2}=\delta$ then $\chi$ is a finite
  morphism~\cite[Chapter I, \S 5, Theorem~7]{sh77} and onto of
  degree $\deg\cLG (n)$~\cite[Corollary~5.6]{mu76} Similarly for
  $\chi'$.
\end{proof}

The set $\cZ\cap \cLG (n)$ is sometimes referred to as the {\em
  base locus}. The interesting part of the theorem occurs when
the base locus $\cZ \cap \cLG (n)=\emptyset$ since in this
situation very specific information on the number of solutions is
provided.  If $\cZ \cap \cLG (n)=\emptyset$ and
$\binom{n+1}{2}=\delta$ (or $n(n+1)=\delta$) then one says that
$\chi$ (or $\chi'$) describes a {\em finite morphism} from the
projective variety $\cLG (n)$ onto the projective space
$\pn{\delta}$ (or $\pn{\delta/2})$.

This last situation is most desirable and this motivates the
following definition.

\begin{defn}                                      \label{defi-degen}
  A particular symmetric transfer function~$G(s)$ is called {\em
    nondegenerate} if $\cZ \cap \cLG (n)=\emptyset$. A system
  which is not nondegenerate will be called degenerate.
\end{defn}

In terms of matrices a symmetric transfer function
$G(s)=D(s)^{-1}N(s)$ is degenerate as soon as there is a
Lagrangian subspace $\row [F_1\ F_2]\in \cLG (n)$, such that
$$
\det\vier{D(s)}{N(s)}{F_1}{F_2}=0.
$$

In a slightly more geometric language this means that the
Hermann-Martin curve~\cite{ma78} defined by $\row [D(s)\ N(s)]$
is fully contained in a Lagrangian hyper-plane defined by $\row
[F_1\ F_2]$. In the study of the static pole placement
problem~\cite{br81} and the dynamic pole placement
problem~\cite{ro94} definitions analogous to
Definition~\ref{defi-degen} played an important role.

The next lemmas give specific information under what conditions
$\cZ \cap \cLG (n)=\emptyset$, i.e. under what conditions a
symmetric transfer function is nondegenerate. Similar results
were crucial in proving the pole placement results
in~\cite{br81,ki04,ro94}.
\begin{lem}                           \label{lemma2.5}
  If $\delta< \binom{n+1}{2}=\dim\cLG (n)$ then every $n\times n$
  symmetric transfer function of McMillan degree $\delta$ is
  degenerate. Similarly, any $n\times n$ Hamiltonian transfer
  function of McMillan degree $\delta$ is degenerate, if $\delta<
  n(n+1)=2\dim\cLG (n)$
\end{lem}
\begin{proof}
  $\dim\cZ\geq N-\delta-1$ as $\cZ$ is defined by $\delta+1$
  linear equations ( $\delta/2+1$ many in the Hamiltonian case).
  If $\dim\cLG (n) >\delta$ (or $\dim\cLG (n) >\delta/2$ in the
  Hamiltonian case), then $\cZ \cap \cLG (n)$ is nonempty by the
  (projective) dimension theorem (see e.g.~\cite[Chapter~I,
  Theorem~7.2]{ha77}).
\end{proof}

\begin{lem}                           \label{gen-non}
  If $\delta = \binom{n+1}{2}=\dim\cLG (n)$ (or $\delta =
  n(n+1)$), then a generic set of $n\times n$ symmetric (or
  Hamiltonian) transfer function of McMillan degree $\delta$ is
  nondegenerate.
\end{lem}
\begin{proof}
  Let $\mathcal{Q}$ be the set of all $n\times n$ symmetric
  transfer functions of McMillan degree $n$. $\mathcal{Q}$ can be
  given the structure of a quasi-projective variety. For this
  recall the definition of the projective variety
  $K^\delta_{n,n}$ introduced in~\cite{ro94} and which
  compactifies the set of all $n\times n$ transfer functions of
  McMillan degree $\delta$. An element (Hermann-Martin curve)
  $\row [D(s)\ N(s)]\in K^\delta_{n,n}$ describes an element of
  $\mathcal{Q}$ as soon as $\deg\det D(s)=\delta$ and
  $D(s)N(s)^t=N(s)D(s)^t$. The last condition translates into
  some linear conditions to be satisfied among the Pl\"ucker
  coordinates of $K^\delta_{n,n}$. The resulting sub-variety of
  $K^\delta_{n,n}$ constitutes a natural compactification of
  $\mathcal{Q}$ and $\mathcal{Q}$ itself is a quasi-projective
  subset.
  
  Consider now the coincidence set
  $$
  {\mathcal S}:=\left\{ (D(s)^{-1}N(s);\ 
    F_1,F_2)\in\mathcal{Q}\times \cLG (n) \mid
    \det\vier{D(s)}{N(s)}{F_1}{F_2}=0 \right\} .
  $$
  Since $\cLG (n)$ is projective the projection onto
  $\mathcal{Q}$ is an algebraic set by the main theorem of
  elimination theory (see e.g.~\cite{mu76}). The set of
  nondegenerate systems forms therefore a Zariski open subset of
  $\mathcal{Q}$.  We have shown the result if we can exhibit one
  $n\times n$ transfer function of McMillan degree
  $\binom{n+1}{2}$ which is nondegenerate. The next lemma gives
  such an example and the claim therefore follows. Note that the
  previous arguments run completely similar for the Hamiltonian
  case and it therefore remains to construct one example as well.
  However, the symmetric Hamiltonian transfer function $G(s^2)$
  does exactly the job.

\end{proof}

\begin{lem}     \label{exa-non}
  The symmetric transfer function $G(s):=\left[
  \begin{array}{cccc}
      \frac{1}{s}  & &&\\
      &\frac{1}{s^2}   &&\\
      &&\ddots&\\
      & &&\frac{1}{s^n}
\end{array}
\right]$ is nondegenerate.
\end{lem}

\begin{proof}
  First it is clear that $G(s)$ has McMillan degree
  $\delta=\binom{n+1}{2}$ and that
  $$
  [D(s)\ N(s)]= \left[
  \begin{array}{cccc|cccc}
      s  &   &      &    & 1&  &      &\\
      &s^2&      &    &  & 1&      &\\
      &   &\ddots&    &  &  &\ddots&\\
      & & &s^n & & & &1
    \end{array}
  \right]
  $$
  forms a left coprime factorization of $G(s)$. Let
  $$
  R:=\left[
  \begin{array}{cccc}
      &&&1   \\
      & & .\cdot   &\\
      &  .\cdot   &&\\
      1 &&&
   \end{array}
 \right]
 $$
 and assume by contradiction that $G(s)$ is degenerate. It
 therefore exists $\row [F_1\ F_2]\in \cLG (n)$, such that
 \begin{equation}                         \label{sing-det}
 0=\det\vier{D(s)}{N(s)}{F_1}{F_2}=\det\vier{D(s)}{N(s)R}{F_1}{F_2R}.
 \end{equation}
 Let $S\in Gl_n$ be the matrix which transforms the $n\times 2n$
 matrix $[F_1\ F_2R]$ into row reduced echelon form, i.e.
\begin{equation}
  \label{cell}
 [(SF_1)\ \ (SF_2R)]=
\left[
\begin{array}{cccccccccccccccc}
\ast &\cdots&\ast&1 &0   &\cdots&0 &0
      &\cdots&0   &\cdots&0   &0&0&\cdots&0\\
      \ast &\cdots&\ast&0 &\ast&\cdots&\ast &1
      &\cdots&0   &\cdots&0   &0&0&\cdots&0\\
      \vdots& &\vdots&\vdots&\vdots& &\vdots&\vdots&
      &\vdots&    &\vdots&\vdots&\vdots&&\vdots\\
      \ast &\cdots&\ast&0 &\ast&\cdots&\ast
      &0&\cdots&\ast&\cdots&\ast&1 &0&\cdots&0
  \end{array}
\right]=:[\tilde{F_1}\ \ \tilde{F_2}]
\end{equation}
Let
$$
i_1< \ldots i_k\leq n < i_{k+1} <\ldots i_n \leq 2n
$$
be the pivot indices. We claim that the first $k$ pivot
indices determine the last $n-k$ pivot indices uniquely.  For
this let $\hat{i}_1<\ldots<\hat{i}_{n-k}$ be the complementary
indices of the indices $\{i_1, \ldots i_k\}$ inside the set
$\{1,\ldots,n\}$. Then we claim that:
\begin{eqnarray*}
i_{k+1}&=& 2n-\hat{i}_{n-k}+1\\
    &\vdots&  \\
i_{n}&=& 2n-\hat{i}_{1}+1.
\end{eqnarray*}
Indeed, if this is not the case then it follows that
$\tilde{F_1}R(\tilde{F_2})^t$ cannot be symmetric for any choice
of values in the row reduced echelon form\eqr{cell}. On the other
hand the matrix $\tilde{F_1}R(\tilde{F_2})^t$ has to be symmetric
since by assumption $F_1(F_2)^t$ is symmetric.

The indices $i_1,\ldots,i_n$ describe the maximal Pl\"ucker
coordinate (with regard of the Bruhat order) of $\row [F_1\ \ 
F_2R]$ which is nonzero and the correponding cofactor of $[D(s) \ 
\ N(s)R]$ is computed as $\pm s^\alpha$, where
$\alpha=\sum_{\ell=1}^{n-k}\hat{i}_\ell$. In general there are
other fullsize minors (Pl\"ucker coordintes) of $[D(s) \ \ 
N(s)R]$ which have the form $\pm s^\alpha$. All other Pl\"ucker
coordinates with this value are however not comparable with
regard to the Bruhat order and since $i_1,\ldots,i_n$ was the
maximal nonzero Pl\"ucker coordinate of $\row [F_1\ \ F_2R]$ it
follows that the determinant expansion in\eqr{sing-det} cannot be
zero. This is a contradiction and it follows that $G(s)$ is
nondegenerate.
\end{proof}

\begin{rem}
  For the static pole placement problem Brockett and
  Byrnes~\cite{br81} showed that the osculating normal curve
  $$
  \row \left[
   \begin{array}{cccccc}
    1 & s& s^2&\ldots &\ldots&s^{m+p-1}\\
    0 & 1& 2s &\ldots&\ldots &\binom{m+p-1}{1}s^{m+p-2}\\
    \vdots& & \ddots &&&\vdots\\
    0&\ldots&0 & 1&\ldots & \binom{m+p-1}{m-1}s^{p}
  \end{array}
\right]\in \G(m,\K^{m+p})
$$
is nondegenerate. Also in this situation the Pl\"ucker
coordinates have the simple form $\pm s^\beta$, where
$\beta=\sum_{\ell=1}^{m}i_\ell-\ell$ and there are no two
Pl\"ucker coordinates which are comparable in the Bruhat order
and give rise to the same monomial $s^\beta$.
\end{rem}

We have now all pieces together in order to prove the main
result.

\begin{proof}[Proof of Theorem~\ref{main}]
  W.l.o.g. we focus on the case of symmetric transfer functions.
  The arguments for the Hamiltonian case run completely similar.
  Note however, that the closed loop characteristic polynomial of
  a Hamiltonian system is always an {\em even} polynomial.
  Therefore our definition of generic pole-assignability for
  Hamiltonian systems restricts to the space of even polynomials.
  Since the dimension of the space of even monic polynomials of
  degree $\delta$ is $\delta/2$, the appropriate condition for
  Hamiltonian systems is $\delta/2 \leq \binom{n+1}{2}$. With the
  comments in mind we return to the proof for symmetric transfer
  functions.
  
  When $\delta> \binom{n+1}{2}$ then a simple dimension argument
  shows that the image of the characteristic map $\chi$ described
  in\eqr{central} has dimension at most $\binom{n+1}{2}$ and
  therefore there is a Zariski open set in $\pn{\delta}$ not in
  the image of~$\chi$.

  When $\delta= \binom{n+1}{2}$ then Lemmas~\ref{gen-non}
  and~\ref{exa-non} show that there is a generic set of $n\times
  n$ symmetric transfer functions of McMillan degree $\delta$
  which are nondegenerate. The characteristic map\eqr{central}
  has therefore no base locus and every point in the image of
  $\chi$ has $\deg\cLG (n)$ pre-image points when counted with
  multiplicities. The degree of the variety $\cLG (n)$ was
  recently computed by Totaro~\cite{to03} and it resulted in the
  number\eqr{deg-formula}.
  
  A priori the geometric formulation only predicts $\deg\cLG (n)$
  many solutions inside $\cLG (n)$ and it is not clear if all
  these solutions correspond to regular feedback laws of the form
  $u=-Fy+v$. If $G(s)$ is a strictly proper symmetric transfer
  function then this is indeed the case and the same argument
  applies as in~\cite{br81}.
\end{proof}

\nocite{he89a1,wa96}
\def\cprime{$'$} \def\cprime{$'$}

\end{document}